\documentclass[10pt]{article}
\usepackage{latexsym, amscd, 
amssymb}
    \title{{\bf  Differential equations and conformal field theories}}
    \author{Yi-Zhi Huang}
    \date{}
    \begin{document}
    \bibliographystyle{alpha}
    \maketitle

\newtheorem{thm}{Theorem}[section]
\newtheorem{defn}[thm]{Definition}
\newtheorem{prop}[thm]{Proposition}
\newtheorem{cor}[thm]{Corollary}
\newtheorem{lemma}[thm]{Lemma}
\newtheorem{rema}[thm]{Remark}
\newtheorem{app}[thm]{Application}
\newtheorem{prob}[thm]{Problem}
\newtheorem{conv}[thm]{Convention}
\newtheorem{conj}[thm]{Conjecture}
\newcommand{\halmos}{\rule{1ex}{1.4ex}}
\newcommand{\pfbox}{\hspace*{\fill}\mbox{$\halmos$}}
\newcommand{\binom}[2]{{{#1}\choose {#2}}}
\newcommand{\text}[1]{\mbox{\rm #1}}
\newcommand{\mod}{\;\;\mbox{\rm mod}\;}

 \newcommand{\nn}{\nonumber \\}
 
	\newcommand{\nno}{\nonumber}
	\newcommand{\lbar}{\bigg\vert}
\newcommand{\mbar}{\mbox{\large $\vert$}}
	\newcommand{\p}{\partial}
	\newcommand{\dps}{\displaystyle}
	\newcommand{\bra}{\langle}
	\newcommand{\ket}{\rangle}
 \newcommand{\res}{\mbox{\rm Res}}
\renewcommand{\hom}{\mbox{\rm Hom}}
\newcommand{\hol}{\mbox{\rm Hol}}
\newcommand{\dt}{\mbox{\rm Det}}
\newcommand{\edo}{\mbox{\rm End}\;}
 \newcommand{\pf}{{\it Proof.}\hspace{2ex}}
 \newcommand{\epf}{\hspace*{\fill}\mbox{$\halmos$}}
 \newcommand{\epfv}{\hspace*{\fill}\mbox{$\halmos$}\vspace{1em}}
 \newcommand{\epfe}{\hspace{2em}\halmos}
\newcommand{\nord}{\mbox{\scriptsize ${\circ\atop\circ}$}}
\newcommand{\wt}{\mbox{\rm wt}\ }
\newcommand{\tr}{\mbox{\rm Tr}}
\newcommand{\swt}{\mbox{\rm {\scriptsize wt}}\ }
\newcommand{\clr}{\mbox{\rm clr}\ }

\begin{abstract}
We discuss the recent results of the author 
on the existence of systems of differential equations for 
chiral genus-zero and 
genus-one correlation functions in conformal field theories. 
\end{abstract}

\renewcommand{\theequation}{\thesection.\arabic{equation}}
\renewcommand{\thethm}{\thesection.\arabic{thm}}
\setcounter{equation}{0}
\setcounter{thm}{0}

\section{Introduction}

Two-dimensional conformal field theories form a particular class of
nontopological quantum field theories which have now been formulated and
studied rigorously using various methods from different branches of
mathematics. In physics, these theories describe perturbative string
theory and also critical phenomena in condensed matter physics.  They
are also used to describe such phenomena as disorder in condensed matter
physics and to construct nonperturbative objects such as $D$-branes in
string theory.  In mathematics, they are closely related to
infinite-dimensional Lie algebras, infinite-dimensional integrable
systems, the Monster (the largest finite sporadic simple group), modular
functions and modular forms, Riemann surfaces and algebraic curves, knot
and three-manifold invariants, Calabi-Yau manifolds and mirror symmetry,
and many other branches of mathematics. We also expect that many
mathematical problems can be solved by constructing and studying the
corresponding conformal field theories. Moreover, 
the study of such theories might provide hints to possible
deep connections among these different branches of mathematics
and will probably shed light on the construction and study of 
higher-dimensional nontopological quantum field theories.

Mathematically, a geometric formulation of conformal field theory was
first given around 1987 by Segal \cite{S1} \cite{S2} \cite{S3} and
Kontsevich.  In \cite{S2} and \cite{S3}, Segal further introduced the
important notions of modular functor and weakly conformal field theory
which describe mathematically the subtle and deep chiral structures
in conformal field theories.
One urgent problem is to give a
construction of (chiral) conformal field theories in this sense. 

To construct conformal field theories in this sense and to study these
conformal field theories, it is necessary to construct and study chiral
correlation functions on Riemann surfaces. For chiral correlation
functions on genus-zero Riemann surfaces (or simply called chiral
genus-zero correlation functions) associated to lowest weight vectors in
minimal models \cite{BPZ} and in Wess-Zumino-Novikov-Witten models
\cite{W}, Belavin-Polyakov-Zamolodchikov and Knizhnik-Zamolodchikov
found in their seminal works \cite{BPZ} and \cite{KZ}, respectively,
that these functions actually satisfy certain systems of differential
equations of regular singular points (now called the BPZ equations and
the KZ equations, respectively). In the case of
Wess-Zumino-Novikov-Witten models, it is also known from the works of
Tsuchiya-Ueno-Yamada \cite{TUY} and Bernard \cite{B1} \cite{B2} that
chiral correlation functions on higher-genus Riemann surfaces (or simply
called chiral higher-genus correlation functions) satisfy systems of
differential equations of KZ type.  These equations play fundamental
roles in the construction and study of the minimal models and
Wess-Zumino-Novikov-Witten models. 

A natural question is whether for general conformal field theories
satisfying natural conditions, there exist systems of differential
equations of regular singular points satisfied by chiral genus-zero
correlation functions. More generally, we are interested in whether
there exist systems of differential equations for 
chiral higher-genus
correlation functions.  The existence of such equations will allow us to
study chiral correlation functions using the theory of differential
equations and to construct chiral conformal field theories using these
correlation functions.

Recently, in \cite{H8} and \cite{H10}, the 
author established the existence of such differential
equations in the genus-zero and genus-one cases under suitable natural
conditions and applied these equations to the construction of genus-zero
and genus-one chiral theories.  In the present paper, after a brief
discussion of the notion of conformal field theories in the sense of
Segal and Kontsevich, we give an overview of these differential
equations. For details, see \cite{H8} and \cite{H10}.  For a recent
exposition on conformal field theories in the sense of Segal and
Kontsevich and the author's program of constructing such theories from
representations of vertex operator algebras, see \cite{H9}. 

In the next section, we recall roughly 
what a conformal field theory is in the sense of Segal 
\cite{S1} \cite{S2} \cite{S3} and
Kontsevich and what a weakly conformal field theory is in the sense of 
Segal \cite{S2} \cite{S3}. We discuss systems of differential equations for 
chiral genus-zero and genus-one correlation functions in Sections 3 and 4,
respectively.

\paragraph{Acknowledgment}  This research 
is supported in part 
by NSF grant DMS-0070800. I am grateful to Sen Hu and Jing Song He for their 
hospitality during the conference.

\renewcommand{\theequation}{\thesection.\arabic{equation}}
\renewcommand{\thethm}{\thesection.\arabic{thm}}
\setcounter{equation}{0}
\setcounter{thm}{0}

\section{Conformal field theories}

Consider the following geometric category:
The objects are  ordered finite sets of copies of $S^{1}$.
The morphisms are  conformal 
equivalence classes of Riemann surfaces whose boundary
components are analytically parametrized by the copies of $S^{1}$ in their
domains and codomains.  The compositions of 
morphisms are given using the boundary parametrizations in the obvious
way. This category has a symmetric monoidal category structure 
for which the monoidal structure is defined
by disjoint unions of objects and morphisms.

Roughly speaking, a {\it conformal field theory} is a projective linear 
representation of this category, that is, a locally 
convex topological vector space $H$ (called the {\it state space})
with a nondegenerate bilinear
form and a projective functor from this category to the symmetric monoidal 
category with traces generated by $H$ (that is, the category whose 
objects are tensor powers of $H$ and  morphisms are
trace-class maps), satisfying some natural conditions. 

Conformal field theories in
general have holomorphic (or chiral) and antiholomorphic (or antichiral)
parts.  Both parts also
satisfy an axiom system which defines {\it weakly conformal field
theories}. Roughly speaking, weakly conformal field
theories are representations of geometric
categories obtained from holomorphic vector bundles over the moduli space
of Riemann surfaces with parametrized boundaries.

Our strategy is to construct chiral or antichiral 
genus-zero and genus-one parts of conformal field theories first and
then using these to construct the full conformal field theories.  In the
remaining part of this paper, we shall discuss only chiral genus-zero and 
genus-one theories.

\renewcommand{\theequation}{\thesection.\arabic{equation}}
\renewcommand{\thethm}{\thesection.\arabic{thm}}
\setcounter{equation}{0}
\setcounter{thm}{0}

\section{Differential equations and chiral genus-zero correlation functions}

We first explain the main ingredients of chiral or antichiral genus-zero
theories.  The chiral or antichiral parts of genus-zero theories have
been shown to be essentially equivalent to algebras of intertwining
operators among modules for suitable vertex operator algebras
(see \cite{H2}--\cite{H5} and \cite{H7}). So here
we briefly describe vertex operator algebras, modules and intertwining
operators. 

A vertex operator algebra is 
a $\mathbb{Z}$-graded vector space $V=\coprod_{n\in \mathbb{Z}}V_{(n)}$
equipped with
a {\it vertex operator map} $Y:V\otimes V\to V[[z, z^{-1}]]$, 
a {\it vacuum} $\mathbf{1}\in V$ and a {\it conformal element} $\omega\in V$,
satisfying a number of axioms. One version of the main axiom
is the following: For $u_{1}, u_{2}, v\in V$, 
$v'\in V'=\coprod_{n\in \mathbb{Z}}V_{(n)}^{*}$, the series 
\begin{eqnarray*}
&\langle v', Y(u_{1}, z_{1})Y(u_{2}, z_{2})v\rangle&\\
&\langle v', Y(u_{2}, z_{2})Y(u_{1}, z_{1})v\rangle&\\
&\langle v', Y(Y(u_{1}, z_{1}-z_{2})u_{2}, z_{2})v\rangle&
\end{eqnarray*}
are absolutely convergent
in the regions $|z_{1}|>|z_{2}|>0$, $|z_{2}|>|z_{1}|>0$ and
$|z_{2}|>|z_{1}-z_{2}|>0$, respectively, to a common rational function 
in $z_{1}$ and $z_{2}$ with the only possible poles at $z_{1}, z_{2}=0$ and 
$z_{1}=z_{2}$. 
Other axioms are:
$$\dim V_{(n)}<\infty,$$
$$V_{(n)}=0$$
when $n$ is sufficiently negative (these are called 
the {\it grading-restriction conditions});
for $u, v\in V$, $Y(u, z)v$ contain only finitely many negative 
power terms; for $u\in V$, 
$$Y(\mathbf{1}, z)=1,\;\;\;\;\lim_{z\to 0}Y(u, z)\mathbf{1}=u;$$
let $L(n): V\to V$ be defined by 
$Y(\omega, z)=\sum_{n\in \mathbb{Z}}L(n)z^{-n-2}$,
then 
$$[L(m), L(n)]=(m-n)L(m+n)+\frac{c}{12}(m^{3}-m)\delta_{m+b, 0}$$
($c$ is called the {\it central charge} of $V$),
$$\frac{d}{dz}Y(u, z)=Y(L(-1)u, z)
\;\;\;\;\;\;\mbox{\rm for}\;\;u\in V$$
and
$$L(0)u=nu \;\;\;\;\;\;\mbox{\rm for}\;\;u\in V_{(n)}$$
($n$ is called the {\it weight} of $u$ and is denoted $\wt u$).
For $u\in V$, we write $Y(u, z)=\sum_{n\in \mathbb{Z}}u_{n}z^{-n-1}$
where $u_{n}\in \edo\; V$.

A {\it $V$-module} is an $\mathbb{C}$-graded vector space 
$W=\coprod_{n\in \mathbb{C}}W_{(n)}$ equipped with a vertex operator map 
$Y_{W}: V\otimes W\to W[[z, z^{-1}]]$ satisfying all those axioms for $V$
which still make sense. Let $W_{1}$, $W_{2}$ and $W_{3}$ be 
$V$-modules. 
An {\it intertwining operators} of type $\binom{W_{3}}{W_{1}W_{2}}$
is a linear map $\mathcal{Y}: W_{1}\otimes
W_{2}\to W_{3}\{z\}$, where 
$W_{3}\{z\}$ is the space of all series in complex powers of $z$ with
coefficients in $W_{3}$,  satisfying all those axioms 
for $V$ which still 
make sense. That is, for $w_{1}\in W_{1}$ and $w_{2}\in W_{2}$,
the real parts of the 
powers of $z$ in nonzero terms 
in the series $\mathcal{Y}(w_{1}, z_{2})w_{2}$ 
have a lower bound;
for $u\in V$, $w_{1}\in W_{1}$, $w_{2}\in W_{2}$
and $w'_{3}\in W'_{3}=\coprod_{n\in \mathbb{C}}(W_{3})_{(n)}^{*}$,
\begin{eqnarray*}
&\langle w'_{3}, Y_{W_{3}}(u, z_{1})\mathcal{Y}(w_{1}, z_{2})w_{2}\rangle&\\
&\langle w'_{3}, \mathcal{Y}(w_{1}, z_{2})Y_{W_{2}}(u, z_{1})w_{2}\rangle&\\
&\langle w'_{3}, \mathcal{Y}(Y_{W_{1}}(u, z_{1}-z_{2})w_{1}, z_{2})
w_{2}\rangle&
\end{eqnarray*}
are absolutely convergent
in the regions $|z_{1}|>|z_{2}|>0$, $|z_{2}|>|z_{1}|>0$ and
$|z_{2}|>|z_{1}-z_{2}|>0$, respectively, to a common (multivalued)
analytic function 
in $z_{1}$ and $z_{2}$ with the only possible singularities
(branch points)
at $z_{1}, z_{2}=0$ and 
$z_{1}=z_{2}$; also
$$\frac{d}{dz}\mathcal{Y}(w_{1}, z)=Y(L(-1)w_{1}, z).$$
For more details on basic notions and properties in the theory of
vertex operator algebras, see \cite{FLM} and \cite{FHL}.

We need the following notions to state the result on differential equations 
in the 
genus-zero case: Let $V$ be a vertex 
operator algebra and $W$ a $V$-module. Let $C_{1}(W)$
be the subspace of $W$ spanned by elements of the form $u_{-1}w$ for
$u\in V_{+}=\coprod_{n>0}V_{(n)}$ and $w\in W$. If $\dim
W/C_{1}(W)<\infty$, we say that $W$ is {\it $C_{1}$-cofinite} or $W$
satisfies the {\it $C_{1}$-cofiniteness condition}.

For chiral genus-zero theories, the main objects
we want to construct and study are chiral genus-zero correlation
functions. Let $W_{i}$ for $i=0, \dots, n+1$ and $\widetilde{W}_{i}$
for $i=1, \dots, n-1$ be $V$-modules and let
$\mathcal{Y}_{1}, \mathcal{Y}_{2},
\dots, \mathcal{Y}_{n-1}, \mathcal{Y}_{n}$ be intertwining operators
of types 
${W'_{0}\choose W_{1}\widetilde{W}_{1}}$, ${\widetilde{W}_{1}\choose 
W_{2}\widetilde{W}_{2}}, \dots, 
{\widetilde{W}_{n-2}\choose W_{n-1}\widetilde{W}_{n-1}}$, 
${\widetilde{W}_{n-1}\choose W_{n}W_{n+1}}$, 
respectively.
Let $w_{i}\in W_{i}$ for $i=0, \dots, n+1$.
Formally, chiral genus-zero correlation
functions are given by 
series of the form
\begin{equation}\label{p-prod}
\langle w_{0}, \mathcal{Y}_{1}(w_{1}, z_{1})\cdots
\mathcal{Y}_{n}(w_{n}, z_{n})w_{n+1}\rangle.
\end{equation}

\begin{thm}\label{sys2}
Let $W_{i}$ for $i=0, \dots, n+1$ be $V$-modules satisfying the 
$C_{1}$-cofiniteness condition. Then 
for any $w_{i}\in W_{i}$ for $i=0, \dots, n+1$, there exist
$$a_{k, \;l}(z_{1}, \dots,  z_{n})\in
\mathbb{C}[z_{1}^{\pm 1}, \dots, z_{n}^{\pm 1},
(z_{1}-z_{2})^{-1}, (z_{1}-z_{3})^{-1}, \dots, (z_{n-1}-z_{n})^{-1}],$$
for $k=1, \dots, m$ and
$l=1, \dots, n,$
such that for any $V$-modules $\widetilde{W}_{i}$ for $i=1, \dots, n-1$, any
intertwining operators $\mathcal{Y}_{1}, \mathcal{Y}_{2},
\dots, \mathcal{Y}_{n-1}, \mathcal{Y}_{n}$,
of types 
${W'_{0}\choose W_{1}\widetilde{W}_{1}}$, ${\widetilde{W}_{1}\choose 
W_{2}\widetilde{W}_{2}}, \dots, 
{\widetilde{W}_{n-2}\choose W_{n-1}\widetilde{W}_{n-1}}$, 
${\widetilde{W}_{n-1}\choose W_{n}W_{n+1}}$, 
respectively, the series (\ref{p-prod})
satisfy the expansions of the system 
of differential equations
$$\frac{\partial^{m}\varphi}{\partial z_{l}^{m}}+
\sum_{k=1}^{m}a_{k, \;l}(z_{1}, \dots,  z_{n})
\frac{\partial^{m-k}\varphi}{\partial z_{l}^{m-k}}=0,\;\;\;l=1, \dots, n$$
in the region $|z_{1}|>\cdots |z_{n}|>0$. Moreover, there exist
$a_{k, \;l}(z_{1}, \dots,  z_{n})$ for $k=1, \dots, m$ and
$l=1, \dots, n$ such that the singular points of the 
corresponding system are regular.
\end{thm}

Similar systems of differential equations have also been obtained by
Nagatomo and Tsuchiya in \cite{NT}. 

Using these equations and other results on vertex operator algebras,
modules and intertwining operators,  chiral genus-zero weakly
conformal field theories in the sense of Segal have been constructed. 
In particular,
the direct sum of a complete set of inequivalent irreducible
modules for a suitable vertex operator algebra has a natural structure
of an intertwining operator algebra. Thus for such 
a vertex operator algebra, (\ref{p-prod}) is
absolutely convergent in the region $|z_{1}|>\cdots |z_{n}|>0$ and 
associativity and commutativity for intertwining operators hold.  For
more details on intertwining operator algebras and chiral genus-zero
weakly conformal field theories, see \cite{H1}--\cite{H9}.

\renewcommand{\theequation}{\thesection.\arabic{equation}}
\renewcommand{\thethm}{\thesection.\arabic{thm}}
\setcounter{equation}{0}
\setcounter{thm}{0}

\section{Differential equations and chiral genus-one correlation functions}

The second logical step is to construct chiral genus-one theories, that
is, to construct maps associated to genus-one surfaces and prove the
axioms which make sense for genus-one surfaces.  Assume we have a weakly 
conformal 
field theory in the sense of Segal. Then for given elements in the
state space of the theory, these maps give certain functions on the
moduli space of genus-one surfaces with punctures (the space of
conformal equivalence classes of such surfaces).  They can be viewed as
(multivalued)
functions of $z_{1}, \dots, z_{n}\in \mathbb{C}$ and $\tau\in
\mathbb{H}$ (the upper half plane).  Here as usual, $\tau$ corresponds
to a torus given by the parallelogram with vertices $0, 1, \tau$ and
$1+\tau$ and $z_{1}, \dots, z_{n}$ correspond to points on the torus. 

But functions of  $z_{1}, \dots, z_{n}$ and $\tau$ 
are in general only functions on the Teichm\"{u}ller space, not
functions on the moduli space. To construct genus-one theories, we do
need to construct mathematical objects (vector bundles and holomorphic
sections) on the moduli space, not the Teichm\"{u}ller space. The moduli
space is the quotient of the Teichm\"{u}ller space $\mathbb{H}$ by the
modular group $SL(2, \mathbb{Z})$. So we have to construct $SL(2,
\mathbb{Z})$-invariant spaces of functions of the form above.  These
functions in the $SL(2,
\mathbb{Z})$-invariant spaces are called {\it chiral 
genus-one correlation functions}.

The first result in this step was obtained by Zhu \cite{Z}. He 
constructed chiral genus-one correlation functions associated to 
elements of a suitable vertex operator algebra $V$. Using his method,
Miyamoto \cite{M} 
constructed chiral genus-one correlation functions associated to
elements of $V$-modules
among which at most one is not isomorphic to $V$.
But Zhu's method cannot be generalized
to construct chiral genus-one correlation functions associated to  
elements of $V$-modules among which at least two are not isomorphic to $V$,
because he used a recurrence formula 
which cannot be generalized to this general case. 

In \cite{H10}, the author solved completely the problem of constructing 
chiral genus-one correlation functions from chiral genus-zero correlation 
functions. 
As in the genus-zero case, one of the main tools is systems of 
differential equations. 

To construct
these functions from representaions of a vertex operator
algebra, we need some conditions on the vertex operator algebra
and its modules.
We first need some concepts:
Let $V$ be a vertex operator algebra and $W$ a $V$-module.
Let $C_{2}(W)$ be the subspace of $W$ spanned by elements of the form 
$u_{-2}w$ for $u, w\in W$.  Then we say that $W$ is {\it $C_{2}$-cofinite}
or satisfies the {\it $C_{2}$-cofiniteness
condition} if $\dim W/C_{2}(W)<\infty$.  It is easy to see that 
if $V_{(n)}=0$ when $n<0$ and $V_{(0)}=\mathbb{C}\mathbf{1}$,
a $V$-module $W$ is $C_{2}$-cofinite implies $W$ is $C_{1}$-cofinite.

We now assume that the vertex operator algebra $V$ satisfies the condition
that (\ref{p-prod}) is absolutely convergent in the region 
$|z_{1}|>\cdots |z_{n}|>0$ and  associativity and commutativity for 
intertwining operators hold.
We shall use the notation $q_{z}=e^{2\pi i z}$ for $z\in \mathbb{C}$.
Let $W_{i}$, $\widetilde{W}_{i}$ be $V$-modules, and $w_{i}\in W_{i}$ 
for $i=1, \dots, n$.  For any intertwining operators
$\mathcal{Y}_{i}$, $i=1, \dots, n$, of types ${\widetilde{W}_{i-1}\choose
W_{i}\widetilde{W}_{i}}$, respectively, let 
\begin{eqnarray}\label{correl-fn}
\lefteqn{F_{\mathcal{Y}_{1}, \dots, \mathcal{Y}_{n}}(w_{1}, \dots, w_{n}; 
z_{1}, \dots, z_{n}; \tau)}\nn
&&=\tr_{\widetilde{W}_{n}}\mathcal{Y}_{1}(\mathcal{U}(q_{z_{1}})
w_{1}, q_{z_{1}})\cdots
\mathcal{Y}_{n}(\mathcal{U}(q_{z_{n}}) w_{n},
q_{z_{n}})q_{\tau}^{L(0)-\frac{c}{24}},
\end{eqnarray}
where $c$ is the central charge of $V$ and for 
$z\in \mathbb{C}\setminus \{0\}$,
$$U(q_{z})=e^{(\log |q_{z}|+i\arg q_{z})L(0)} 
e^{-L^{+}(A)},\;\;\;0\le \arg q_{z}<2\pi,$$
$L^{+}(A)=\sum_{j\ge 1}A_{j}L(j)$
and $A_{j}$ for $j\ge 1$ are determined by 
\begin{eqnarray*}
\frac{1}{2\pi i}\log(1+2\pi i w)=\left(\exp\left(\sum_{j\in \mathbb{Z}_{+}}
A_{j}w^{j+1}\frac{\partial}{\partial w}\right)\right)w
\end{eqnarray*}

Let
\begin{eqnarray*}
\wp_{1}(z; \tau)&=&\frac{1}{z}+\sum_{(k, l)\ne (0, 0)}
\left(\frac{1}{z-(k\tau+l)}+\frac{1}{k\tau+l}
+\frac{z}{(k\tau+l)^{2}}\right),\\
\wp_{2}(z; \tau)&=&\frac{1}{z^{2}}+\sum_{(k, l)\ne (0, 0)}
\left(\frac{1}{(z-(k\tau+l))^{2}}-\frac{1}{(k\tau+l)^{2}}\right)
=-\frac{\partial}{\partial z}\wp_{1}(z; \tau)
\end{eqnarray*}
be the Weierstrass zeta function and the 
Weierstrass $\wp$-function, respectively, and let $\wp_{m}(z;\tau)$ for $m>2$ 
be the elliptic functions defined recursively by
$$\wp_{m+1}(z, \tau)=-\frac{1}{m}\frac{\partial}{\partial z} \wp_{m}(z;
\tau).$$
We also need the Eisenstein series
\begin{eqnarray*}
G_{2}(\tau)&=&\frac{\pi^{2}}{3}
+\sum_{m\in \mathbb{Z}\setminus \{0\}}\sum_{l\in \mathbb{Z}}
\frac{1}{(m\tau+l)^{2}},\\
G_{2k+2}(\tau)&=&\sum_{(m, l)\ne (0, 0)}
\frac{1}{(m\tau+l)^{2k+2}},\;\;\;\;\;\;\;\;\;\;\;\;k\ge 1.
\end{eqnarray*}
See, for example, \cite{K} and \cite{L}, 
for detailed discussions on 
these elliptic functions and the Eisenstein series. Let
$$R=\mathbb{C}[G_{4}(\tau), G_{6}(\tau),
\wp_{2}(z_{i}-z_{j}; \tau), 
\wp_{3}(z_{i}-z_{j}; \tau)]_{i, j=1, \dots, n,\;i<j},$$
that is, the commutative associative 
algebra over $\mathbb{C}$ generated by the series
$G_{4}(\tau)$, $G_{6}(\tau)$, 
$\wp_{2}(z_{i}-z_{j}; \tau)$ and $\wp_{3}(z_{i}-z_{j}; \tau)$
for $i, j=1, \dots, n$ satisfying $i<j$.
For $m\ge 0$, let $R_{m}$ be the 
subspace of $R$ spanned by elements of the form 
$$G^{k_{1}}_{4}(\tau)G^{k_{2}}_{6}(\tau)
\wp^{k_{3}}_{2}(z_{i}-z_{j}; \tau)
\wp^{k_{4}}_{3}(z_{i}-z_{j}; \tau)$$
for $k_{1}, k_{2}, k_{3}, k_{4}\ge 0$ satisfying 
$4k_{1}+6k_{2}+2k_{3}+3k_{4}=m$. 

We introduce, for any $\alpha\in \mathbb{C}$,
the notation
$$\mathcal{O}_{j}(\alpha)=2\pi i\frac{\partial}{\partial \tau}
+G_{2}(\tau)\alpha
+G_{2}(\tau)\sum_{i=1}^{n}z_{i}\frac{\partial}{\partial z_{i}}
-\sum_{i\ne j}\wp_{1}(z_{i}-z_{j}; \tau)
\frac{\partial}{\partial z_{i}}$$
for $j=1, \dots, n$. We shall also use the notation
$\prod_{j=1}^{m}\mathcal{O}(\alpha_{j})$
to denote the ordered product
$\mathcal{O}(\alpha_{1})\cdots \mathcal{O}(\alpha_{m})$.

Then we have the following result:

\begin{thm}
Let $V$ be a vertex operator algebra satisfying the conditions stated
above and let
$W_{i}$ for $i=1, \dots, n$ be $V$-modules satisfying the 
$C_{2}$-cofiniteness condition. Then 
for any homogeneous $w_{i}\in W_{i}$ ($i=1, \dots, n$), there exist
$$a_{p, \;i}(z_{1}, \dots, z_{n}; \tau)\in 
R_{p},\;\;\; b_{p, \;i}(z_{1}, \dots, z_{n}; \tau)\in R_{2p}$$
for $p=1, \dots, m$ and $i=1, \dots, n$ such that 
for any $V$-modules $\widetilde{W}_{i}$ ($i=1, \dots, n$) and
intertwining operators
$\mathcal{Y}_{i}$   of 
types $\binom{\widetilde{W}_{i-1}}{W_{i}\widetilde{W}_{i}}$ ($i=1, \dots, n$,
$\widetilde{W}_{0}=\widetilde{W}_{n}$), 
respectively, the series 
(\ref{correl-fn})
satisfies the expansion of the system of differential equations
\begin{eqnarray}
&{\displaystyle \frac{\partial^{m}\varphi}{\partial z_{i}^{m}}+
\sum_{p=1}^{m}a_{p, \; i}(z_{1}, 
\dots, z_{n};\tau)
\frac{\partial^{m-p}\varphi}{\partial z_{i}^{m-p}}
=0,}\label{eqn1}\\
&{\displaystyle \prod_{k=1}^{m}\mathcal{O}_{i}\left(\sum_{i=1}^{n}
\wt w_{i}+2(m-k)\right)\varphi\quad\quad\quad\quad\quad\quad\quad
\quad\quad\quad\quad\quad\quad\quad\quad\quad}\nn
&{\displaystyle  
+\sum_{p=1}^{m}
b_{p, \;i}(z_{1}, \dots, z_{n}; \tau)
\prod_{k=1}^{m-p}\mathcal{O}_{i}\left(\sum_{i=1}^{n}
\wt w_{i}+2(m-p-k)\right)
\varphi=0,}\label{eqn2}
\end{eqnarray}
$i=1, \dots, n$, 
in the regions $1>|q_{ z_{1}}|>\cdots >|q_{z_{n}}|>|q_{\tau}|>0$.
Moreover, for fixed $\tau\in \mathbb{H}$,  the singular points of
the (reduced) system (\ref{eqn1}) are regular.
\end{thm}

The elliptic functions and the Eisenstein series 
discussed above have the following modular transformation formulas
(see, for example, 
\cite{K}):
For any 
$$\left(\begin{array}{cc}
\alpha&\beta\\
\gamma&\delta
\end{array}\right)\in SL(2, \mathbb{Z}),$$
\begin{eqnarray*}
G_{2}\left(\frac{\alpha\tau+\beta}{\gamma\tau +\delta}\right)
&=&(\gamma\tau+\delta)^{2}G_{2}(\tau)
-2\pi i\gamma(\gamma\tau+\delta),\label{mod-g2}\\
G_{2k}\left(\frac{\alpha\tau+\beta}{\gamma\tau+\delta}\right)&=&
(\gamma\tau+\delta)^{2k}G_{2k}(\tau),\label{mod-g2k}\\
\wp_{m}\left(\frac{z}{\gamma\tau+\delta}; 
\frac{\alpha\tau+\beta}{\gamma\tau+\delta}\right)&=&
(\gamma\tau+\delta)^{m}\wp_{m}(z, \tau),\label{mod-wpm}
\end{eqnarray*}
for $k\ge 2$ and $m\ge 1$. 
Using these formulas, it is straightforward to verify the following 
modular invariance of the system 
(\ref{eqn1})--(\ref{eqn2}):

\begin{prop}\label{eqn-mod-inv}
Let $\varphi(z_{1}, \dots, z_{n}; \tau)$ be a solution of 
the system (\ref{eqn1})--(\ref{eqn2}). 
Then for any 
$$\left(\begin{array}{cc}
\alpha&\beta\\
\gamma&\delta
\end{array}\right)\in SL(2, \mathbb{Z}),$$
$$\left(\frac{1}{\gamma\tau +\delta}\right)^{\swt w_{1}+\cdots \swt w_{n}}
\varphi\left(\frac{z_{1}}{c\tau +d}, \dots, \frac{z_{n}}{c\tau +d}; 
\frac{\alpha\tau+\beta}{\gamma\tau +\delta}\right)$$
is also a solution of the system (\ref{eqn1})--(\ref{eqn2}).
\end{prop}

Using these systems of differential equations, 
chiral genus-one correlation functions have been constructed as 
the analytic extensions of sums of series of the form (\ref{correl-fn})
in the region $1>|q_{ z_{1}}|>\cdots >|q_{z_{n}}|>|q_{\tau}|>0$.
Together with 
other results in the representation theory of vertex
operator algebras, it has been proved 
that for suitable vertex operator algebras,
the vector space 
of these chiral genus-one correlation functions are invariant under 
a suitable action of the modular group $SL(2, \mathbb{Z})$.
See \cite{H10} for details.

\noindent {\small \sc Department of Mathematics, Rutgers University,
110 Frelinghuysen Rd., Piscataway, NJ 08854-8019 }

\noindent {\em E-mail address}: yzhuang@math.rutgers.edu

\end{document}